\documentclass[hidelinks]{article}
\usepackage{graphicx} 
\usepackage{hyperref}
\usepackage{amsthm}
\usepackage{geometry}
\usepackage{amssymb}
\usepackage{amsmath}
\usepackage{enumitem}
\usepackage{tikz, float}
\usepackage{xcolor}
\definecolor{purple}{rgb}{0.5,0,0.5}
\usepackage{lmodern}
\usepackage[english]{babel}
\usetikzlibrary{arrows.meta,backgrounds}
\definecolor{ink}{HTML}{263A56}
\definecolor{blue}{HTML}{4276B5}
\definecolor{teal}{HTML}{218C88}
\definecolor{rule}{HTML}{C5CFDB}

\newtheorem{theorem}{Theorem}[section]
\newtheorem{lemma}[theorem]{Lemma}
\newtheorem{cor}[theorem]{Corollary}

\theoremstyle{definition}

\newcommand{\claimproofstart}[1][Proof]{\begin{proof}[#1]
\renewcommand{\qedsymbol}{$\boxdot$}}
\newcommand\claimproofend{
\end{proof}
\renewcommand{\qedsymbol}{$\square$}}

\newcommand{\ind}{\mathbf 1}

\newcommand{\cS}{\mathcal{S}}

\newcommand{\cU}{\mathcal{U}}

\title{Positive semidefiniteness of random submatrices}

\author{Rajko Nenadov \footnote{School of Mathematics and Statistics, University of Canterbury, New Zealand. Research supported by the New Zealand Marsden Fund. Email: rajko.nenadov@canterbury.ac.nz}
\and
Lander Verlinde \footnote{School of Mathematics and Statistics, University of Canterbury, New Zealand. \\ Email: lander.verlinde@pg.canterbury.ac.nz}} 
\date{}

\begin{document}

\maketitle

\begin{abstract}
    Given a symmetric matrix $A \in [-1,1]^{n \times n}$, we estimate the probability that a random principal submatrix, obtained by taking each index independently with probability $p$, is positive semidefinite. As a corollary, we show that $O(\varepsilon^{-1})$ indices chosen uniformly at random suffice for testing positive semidefiniteness with squared Frobenius distance threshold $\varepsilon n^2$, thus nearly matching the lower bound of Bakshi, Chepurko, and Jayaram (``Testing Positive Semi-Definiteness via Random Submatrices'', FOCS 2020). The key ingredient is a dense model lemma for positive semidefinite matrices, proved using the hyperplane separation theorem and a contraction inequality for Rademacher random variables.
\end{abstract}

\section{Introduction}
Given an arbitrary matrix $A \in \mathbb{R}^{n \times n}$ and a set $I \subseteq [n]$, we write $A_I$ for the principal submatrix of $A$ indexed by $I$. A natural and widely studied question is how the spectra of $A$ and $A_I$ compare. 
A seminal result in this line of research is the work of Bourgain and Tzafriri \cite{bourgain1987invertibility} on restricted invertibility. They showed that there exist constants $a, c > 0$ such that given an arbitrary matrix~$A$ with columns having Euclidean norm one, there is a set $I \subseteq [n]$ satisfying 
\[ 
    |I| \; \|A\|_{\operatorname{op}}^2 \geq cn, \quad \text{and} \quad  G_I \succeq a^2  \; \textbf{id}_{I},
\]
where $G = A^\top A$ and $\textbf{id}_I$ is the identity matrix indexed by $I$. This result has generated a substantial body of research (see, e.g., \cite{marcus2022interlacingIII,Naor2017restricted,spielman2012elementary,vershynin2001decomposition}).

The study of principal submatrices also played a central role in the Kadison-Singer problem. Anderson's paving formulation of this problem \cite{anderson1979extensions} asked whether for every $\eta \in (0,1)$, there is an integer $r$ such that every self-adjoint, zero-diagonal complex matrix $A$ admits a partition 
\[ [n] = I_1 \cup \dots \cup I_r \]
such that 
\[ \| A_{I_j} \|_{\operatorname{op}} \leq \eta \| A\|_{\operatorname{op}} \text{ for every }j.  \]
The difficulty here lay in the fact that $r$ has to be independent of $n$. Progress on this conjecture was made by Bourgain and Tzafriri \cite{bourgain1991kadison}, proving $O_\eta(\log n)$-part paving. Their proof is based on an intricate random process that requires good control of $(\mathbb{E} \| A_{I_p} \|_{\operatorname{op}}^q)^{1/q} $ for $q = \Theta(\log n)$ and $I_p$ being the binomial random set for some probability $p$. This was further studied by Tropp \cite{tropp2008norms,Tropp2008paving}. Finally, the Kadison-Singer problem (and the paving conjecture) was settled by Marcus, Spielman and Srivastava \cite{marcus2015interlacingII} using interlacing polynomials. In the language of principal submatrices, their main technical ingredient is the following. Given a matrix $A \succeq 0$ with $\max_{i \in [n]} A_{ii} \leq t$ and any integer $r > 0$, a uniformly chosen random partition $I_1 \cup \dots \cup I_r = [n]$, with positive probability, satisfies
\[ \|A_{I_j} \|_{\operatorname{op}} \leq \left( \sqrt{\frac{\|A\|_{\operatorname{op}}}{r}} + \sqrt{t} \right)^2 \text{ for every } j \in \{1, \ldots, r\}.  \]

We continue the study of the spectrum of submatrices by considering the probability that a random principal submatrix has all eigenvalues non-negative, that is, it is positive semidefinite. 
The following is our main result.

\begin{theorem} \label{thm:main_psd}
    Let $A \in [-1,1]^{n \times n}$ be a symmetric matrix which is not positive semidefinite. For $p \in (0,1]$ and a binomial random set $I \sim [n]_p$ obtained by including each index independently with probability~$p$, we have
    $$
        \Pr\left[ A_I \succeq 0 \right] \le \frac{42\sqrt{n}}{ \min_{B \in \cS_n} \| A - B \|_F\, \sqrt{p}},
    $$
    where $\cS_n$ denotes the family of all symmetric positive semidefinite matrices in $[-1,1]^{n \times n}$. 
\end{theorem}

This theorem  has immediate consequences in property testing \cite{Bhattacharjee2024sublinear,needell2022testing,swartworth2025tight}. Bakshi, Chepurko, and Jayaram \cite{bakshi2020testing} studied how efficiently one can distinguish whether a given symmetric matrix $A \in [-1,1]^{n \times n}$ is positive semidefinite or is far away from being so, in the sense
\[ \min_{B \in \cS_n}\|A-B\|_F^2 \ge \varepsilon n^2. \]
They show that a principal submatrix of order $\Tilde{O}(\varepsilon^{-2})$, chosen uniformly at random, suffices to distinguish the two cases with significant probability. They also prove a $\Tilde{\Omega}(\varepsilon^{-1})$ lower bound. As a direct corollary of Theorem \ref{thm:main_psd}, we provide a nearly-matching upper bound.

\begin{cor}
    There exists a constant $C > 0$ such that the following holds. Let $\varepsilon \in (0, 1]$, and suppose $A \in [-1,1]^{n \times n}$ is a symmetric matrix such that
    $$
        \min_{B \in \cS_n} \| A - B \|_F^2\ \ge \varepsilon n^2.
    $$
    Then, for $I \subseteq \{1, \ldots, n\}$ being a set of size $|I| = \min\{n, \lceil C \varepsilon^{-1} \rceil \}$ chosen uniformly at random, we have
    $$
        \Pr\left[ A_I \succeq 0 \right] \le 1/2.
    $$
\end{cor}

\paragraph{Proof overview.}
The proof of Theorem \ref{thm:main_psd} is built upon ideas of dense models, pioneered by Green and Tao \cite{green2008primes}. Dense models provide a way to replace an object supported on a sparse domain by an object on the full domain, while approximately preserving its behaviour after rescaling.

The main idea of the proof is the following. Consider an arbitrary symmetric $A \in [-1,1]^{n \times n} \setminus \cS_n$, and let $d > 0$ be the Frobenius distance from $A$ to a closest $B \in \cS_n$ (recall that $\cS_n$ denotes the family of all symmetric positive semidefinite matrices in $[-1,1]^{n \times n}$). By the hyperplane separation theorem, there exists a symmetric matrix $Z$ with $\|Z\|_F \leq 1$ such that for every $B \in \cS_n$, we have $\langle A - B, Z \rangle \geq d$. For simplicity, suppose $Z$ has zero-diagonal. The proof strategy is to show that if $A_I \succeq 0$, then there exists $B \in \cS_n$ such that $\langle A - B, Z \rangle$ is small, in expectation. This, in turn, gives us a way to compare the probability of $A_I \succeq 0$ with $d$.

Simplifying the sketch further, assume $p = 2^{-m}$ for some $m \in \mathbb{N}$. To find the above described $B \in \cS_n$, we consider $I$ as being obtained by successive uniform random samples $[n] = I_0 \supseteq I_1  \supseteq I_2 \supseteq \ldots \supseteq I_m = I$. That is, each $I_j$ is sampled uniformly at random conditioned on $I_{j-1}$. The matrix $A_{I_j}$ is, in expectation, supported on $1 / 4$ of the entries which are in the support of $A_{I_{j-1}}$. Therefore, in order to have a fair comparison between $A_{I_{j-1}}$ and $A_{I_j}$, we rescale all the entries in $A_{I_j}$ by 4. Importantly, we expect $\langle A_{I_{j-1}}, Z \rangle$ and $\langle 4A_{I_j}, Z \rangle$ to be close. This motivates the following: Starting with $j = 1$, we successively make steps from $4^{j-1} A_{I_{j-1}}$ to $4^{j} A_{I_j}$. Assuming $A_I \succeq 0$, this eventually brings us to the matrix $2^{2m} A_{I_m} \in 2^{2m} \cS_n$. As each such move causes only a small drift with respect to $\langle \cdot, Z \rangle$, we conclude $\langle A, Z \rangle \approx \langle 2^{2m} A_{I_m}, Z \rangle$. Set $B^{(m)} := A_{I_m}$. Then $B^{(m)} \in \cS_n(I_m)$, where $\cS_n(I_j)$ denotes the matrices in $\cS_n$ supported on indices in $I_j$. The goal now is to move from $B^{(m)} \in \cS_n(I_m)$, through $B^{(j)} \in \cS_n(I_j)$ for $m-1 \ge j \ge 1$, back to $B^{(0)} \in \cS_n$, again making only a small drift with respect to $\langle \cdot, Z \rangle$ in each step. This is where the dense models prove their worth. Given a matrix $B^{(j)} \in \cS_n(I_j)$, we show that one can find a `denser' matrix $B^{(j-1)} \in \cS_n(I_{j-1})$, that is, a matrix supported on $I_{j-1}$ which typically contains roughly twice as many indices as $I_j$, such that it approximates $B^{(j)}$ when appropriately scaled,
$$
    \langle 4 B^{(j)}, Z \rangle \approx \langle B^{(j-1)}, Z \rangle.
$$
The following is a schematic drawing of the proof strategy in the case $m = 2$.

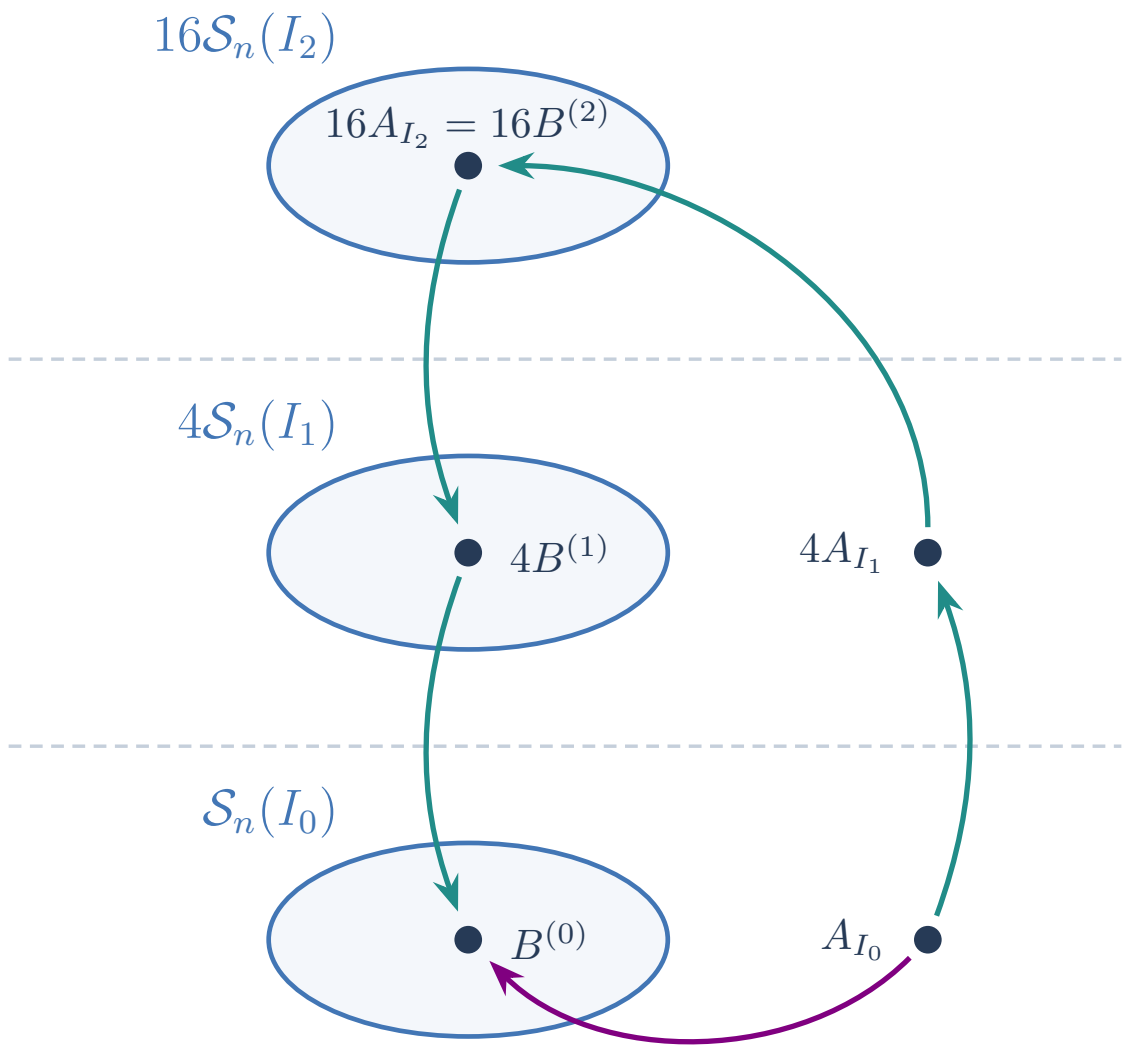
\begin{figure}[H]
\centering
\scalebox{0.8}{%
\begin{tikzpicture}[
    font=\large,
    text=ink,
    region/.style={
        draw=blue,
        fill=blue!6,
        line width=1.1pt
    },
    point/.style={
        circle,
        fill=ink,
        inner sep=2.3pt
    },
    point label/.style={
        font=\normalsize
    },
    arrow/.style={
        -{Stealth[length=3mm,width=2.2mm]},
        draw=teal,
        line width=1.25pt,
        shorten <=6pt,
        shorten >=7pt
    },
    small label/.style={
        font=\small\itshape,
        text=teal,
        inner xsep=2pt,
        inner ysep=3pt,
        outer sep=0pt
    },
    scale label/.style={
        font=\Large,
        text=ink!60
    }
]

\coordinate (B0) at (0,0);
\coordinate (B1) at (0,3.2);
\coordinate (B2) at (0,6.4);
\coordinate (A0) at (3.8,0);
\coordinate (A1) at (3.8,3.2);

\begin{scope}[on background layer]
    \foreach \y in {1.6,4.8} {
        \draw[rule, densely dashed, line width=.8pt]
            (-3.8,\y) -- (5.4,\y);
    }

    \foreach \y in {0,3.2,6.4} {
        \draw[region] (0,\y)
            ellipse[x radius=1.65,y radius=.8];
    }
\end{scope}

\foreach \y/\setlabel in {
    0/{\cS_n(I_0)},
    3.2/{4\cS_n(I_1)},
    6.4/{16\cS_n(I_2)}
} {
    \node[anchor=south east, text=blue]
        at (-.95,{\y+.72}) {$\setlabel$};
}

\draw[arrow]
    (B2) to[out=-110,in=110]
    node[small label, pos=.38, left=2pt] {}
    (B1);

\draw[arrow]
    (B1) to[out=-110,in=110]
    node[small label, pos=.38, left=2pt] {}
    (B0);

\draw[arrow]
    (A0) to[out=70,in=-70]
    node[small label, pos=.62, right=2pt] {}
    (A1);

\draw[arrow]
    (A1) to[out=90,in=0]
    node[small label, pos=.60, above right=2pt] {}
    (B2);

\draw[arrow, draw=purple]
    (A0) to[out=-135,in=-45]
    node[small label, text=purple, midway, below=2pt] {}
    (B0);

\node[point] at (B0) {};
\node[point label, anchor=west]
    at (.22,0) {$B^{(0)}$};

\node[point] at (B1) {};
\node[point label, anchor=west]
    at (.22,3.2) {$4B^{(1)}$};

\node[point] at (B2) {};
\node[point label, anchor=center, inner sep=1pt]
    at (0,6.75) {$16 A_{I_2}= 16B^{(2)}$};

\node[point] at (A0) {};
\node[point label, anchor=east]
    at (3.58,0) {$A_{I_0}$};

\node[point] at (A1) {};
\node[point label, anchor=east]
    at (3.58,3.2) {$4A_{I_1}$};

\end{tikzpicture}%
}
\end{figure}
Following the green arrows, we obtain
\[ |\langle A - B^{(0)} ,Z \rangle | \leq \sum_{j=1}^m 4^{j-1} \; |\langle A_{I_{j-1}} - 4 A_{I_{j}}, Z \rangle| +  \sum_{j=1}^m 4^{j-1} \; |\langle 4 B^{(j)} - B^{(j-1)}, Z \rangle|.  \]
The previous analysis makes use of the assumption that $A_I \succeq 0$. To obtain a bound applicable if this is not the case, we simply take 
\[ d \cdot \ind\{A_I\succeq0\} \leq \sum_{j=1}^{m} 4^{j-1} \Delta_j,  \]
where 
\[  \Delta_j = |\langle A_{I_{j-1}} - 4 A_{I_{j}}, Z \rangle| + L_{Z_{I_{j-1}}}(I_j) \]
and $L_{Z_{I_{j-1}}}(I_j)$ denotes the worst-case guarantee of the dense model lemma (see Lemma \ref{lemma:dense_psd}), that is, the largest (optimal) drift $\langle 4 B^{(j)} - B^{(j-1)}, Z \rangle$ we can get. Hence, to finish the proof, we shall show $\mathbb{E} \Delta_j = O\left(2^{-3j/2} \sqrt{n}\right)$.

\paragraph{Use of AI.} Initial attempts focused on the fingerprint approach to constructing dense models, based on the first author’s recent paper \cite{nenadov2026simplified}. After these attempts proved unsuccessful, ChatGPT Astra suggested a simpler, direct approach using the hyperplane separation theorem and a contraction inequality for Rademacher random variables. The authors have verified and simplified the mathematical claims suggested by the LLM.

\section{Preliminaries}

All matrix inner products are Frobenius inner products, that is, $\langle A, B \rangle = \mathrm{tr}(A^\top B)$. When comparing matrices indexed by different sets, we extend them by zero to the ambient index set. For a matrix $A \in \mathbb{R}^{n \times n}$, we denote by $A_i$ its $i$th column. 
Let $\cU_n$ denote the family of real $n \times n$ matrices with each column having Euclidean norm at most $1$,
$$
    \cU_n = \{A \in \mathbb{R}^{n \times n} \; : \; \|A_i\|_2 \le 1 \text{ for every } i \in [n] \}.
$$
Given $I \subseteq [n]$, set
$$
    \cS_n(I) = \{ A \in \cS_n \; : \; A_i = 0 \text{ for } i \not \in I \},
$$
and define $\cU_n(I)$ analogously. We shall make use of the well-known fact that $\cS_n = \{A^\top A : A \in \cU_n\}$. The same correspondence holds for $\cS_n(I)$ and $\cU_n(I)$. We also note that $\cS_n(I)$ and $\cU_n(I)$ are compact and convex for any fixed $I \subseteq [n]$.

The following is the hyperplane separation theorem we require, and a corollary thereof.

\begin{theorem} \label{thm:hyperplane}
    Let $\mathcal{M} \subseteq \mathbb{R}^{n \times n}$ be a non-empty closed convex set of matrices, and let $A \in \mathbb{R}^{n \times n} \setminus \mathcal{M}$. Let 
    \[d = \min_{B \in \mathcal{M}} \|A - B\|_F > 0.\] 
    Let $B_0 \in \mathcal{M}$ be such that $\|A - B_0 \|_F = d$. Then, $Z = (A - B_0)/d $ satisfies $\langle A-B, Z \rangle \geq d$ for all $B \in \mathcal{M}$.  
\end{theorem}

\begin{lemma} \label{lemma:min_max}
    Let $I \subseteq [n]$ be arbitrary, and consider two arbitrary linear maps $\sigma, \tau: \mathbb{R}^{n \times n} \rightarrow \mathbb{R}^{n \times n}$. Then, for $\cU_n$ and $\cU_n(I)$ as defined above
    \[ \max_{A \in \cU_n(I)} \min_{B \in \cU_n} \| \sigma(A) - \tau(B) \|_F \leq \max_{\substack{U \in \mathbb{R}^{n \times n} \\ \|U\|_F \le 1}} \left( \max_{A \in \cU_n(I)} \langle \sigma(A),U \rangle - \max_{B \in \cU_n} \langle \tau(B), U \rangle \right). \]
\end{lemma}
\begin{proof}
    Consider some $A \in \cU_n(I)$, and let 
    \[ d= \min_{B \in \cU_n} \| \sigma(A) - \tau(B) \|_F. \]
    If $d = 0$, set $U_A$ equal to the all-zero matrix.
    If $d > 0$, observe that $\tau(\cU_n)$ is closed and convex. By Theorem \ref{thm:hyperplane}, there is $U_A \in \mathbb{R}^{n \times n}$ with $\|U_A\|_F \leq 1$ such that 
    \[ \langle \sigma(A) - \tau(B), U_A \rangle \geq d \quad \text{for every } B \in \cU_n. \]
    In particular, 
    \[ \min_{B \in \cU_n} \|\sigma(A) - \tau(B) \|_F \leq \langle \sigma(A), U_A \rangle - \max_{B \in \cU_n} \langle \tau(B),U_A \rangle. \]
    Therefore, 
    \begin{align*}
         \max_{A \in \cU_n(I)} \min_{B \in \cU_n} \| \sigma(A) - \tau(B) \|_F &\leq \max_{A \in \cU_n(I)} \left(  \langle \sigma(A), U_A \rangle - \max_{B \in \cU_n} \langle \tau(B),U_A \rangle \right) \\
         &\leq \max_{\substack{U \in \mathbb{R}^{n \times n} \\ \|U\|_F \le 1}} \left( \max_{A \in \cU_n(I)} \langle \sigma(A),U \rangle - \max_{B \in \cU_n} \langle \tau(B), U \rangle \right) .
    \end{align*}
\end{proof}

The following is a special case of the contraction theorem of Ledoux and Talagrand (see \cite[Theorem 4.12]{ledoux_talagrand}).

\begin{theorem} \label{thm:ledoux}
    Let $A \subseteq \mathbb{R}^n$ be non-empty and bounded, and let $\epsilon = (\epsilon_1, \ldots, \epsilon_n)$ be a vector of independent Rademacher random variables. Writing $|a|=(|a_1|,\ldots,|a_n|)$, we have
    $$
        \mathbb{E} \sup_{a \in A} \bigl| \langle \epsilon, |a| \rangle \bigr| \le 2\, \mathbb{E} \sup_{a \in A} \bigl| \langle \epsilon, a \rangle \bigr|.
    $$
\end{theorem}   

To finish this section, we give a first and easy bound on how a sampled matrix $A_I$ behaves compared to its original matrix.

\begin{lemma} \label{lemma:subset}
    Let $A \in [-1,1]^{n \times n}$ be a zero-diagonal matrix, and let $Z \in \mathbb{R}^{n \times n}$. For $I \subseteq [n]$ chosen uniformly at random, we have
    $$
        \mathbb{E} \left| \langle A, Z \rangle - 4 \langle A_I, Z_I \rangle \right| \le 2 \sqrt{n} \; \| Z \|_F.
    $$
\end{lemma}
\begin{proof}
    For $i \in [n]$, set $\epsilon_i=1$ if $i \in I$ and $\epsilon_i = -1$ otherwise. Note that $\epsilon_i$ are independent Rademacher random variables. For $i \ne j$, set $b_{ij}=A_{ij}Z_{ij}+A_{ji}Z_{ji}$. As $A$ has zero diagonal, we have
    $$
        4\langle A_I,Z_I\rangle-\langle A,Z\rangle
        =\sum_{i<j}b_{ij}(\epsilon_i+\epsilon_j+\epsilon_i\epsilon_j)
        =\sum_i\epsilon_i\sum_{j\ne i}b_{ij}
            +\sum_{i<j}b_{ij}\epsilon_i\epsilon_j.
    $$
    Then
    $$
        \mathbb{E}\left(4\langle A_I,Z_I\rangle-\langle A,Z\rangle\right)^2
        =\sum_i\left(\sum_{j\ne i}b_{ij}\right)^2+\sum_{i<j}b_{ij}^2,
    $$
    as all products of distinct terms have expectation zero. Therefore, by Cauchy-Schwarz,
    \begin{align*}        
        \mathbb{E}\left(4\langle A_I,Z_I\rangle-\langle A,Z\rangle\right)^2 
        &\le (2n-1)\sum_{i<j}b_{ij}^2\\
        &\le (2n-1)\sum_{i\ne j} 2 Z_{ij}^2
        \le 4n\|Z\|_F^2,
    \end{align*}
    where the second inequality uses $|A_{ij}|\le1$. A final application of Cauchy-Schwarz proves the lemma.
\end{proof}

\section{Constructing the dense model}
As mentioned in the introduction, the key technical ingredient is a dense model lemma for symmetric positive semidefinite matrices. This allows us to move back to the unscaled, dense universe and find the matrix $B^{(0)}$. In the following statement, one should think of $X$ as the matrix $B^{(j+1)}$, and $Y$ as the desired $B^{(j)}$.

Following ideas from Bender and Krauthgamer \cite{bender2026optimal}, we prove the dense model lemma using the corollary of the hyperplane separation theorem (Lemma \ref{lemma:min_max}) and the contraction theorem of Ledoux and Talagrand (Theorem \ref{thm:ledoux}).

\begin{lemma} \label{lemma:dense_psd}
    Let $Z \in \mathbb{R}^{n \times n}$ be a symmetric, zero-diagonal matrix. For $I \subseteq [n]$ chosen uniformly at random, we have
    $$
        \mathbb{E} L_{Z}(I) \le 10 \sqrt{n} \; \|Z\|_F,
    $$
    where
    $$
        L_{Z}(I) = \max_{X \in \cS_n(I)} \min_{Y \in \cS_n} \left| \langle 4 X - Y, Z \rangle \right|.
    $$
\end{lemma}
\begin{proof}
    Using $\cS_n = \{A^\top A : A \in \cU_n\}$ and the symmetry of $Z$, we have
    $$
        L_Z(I)
        = \max_{A \in \cU_n(I)} \min_{B \in \cU_n} \left|\langle 4 A^\top A - B^\top B, Z \rangle\right|
        = \max_{A \in \cU_n(I)} \min_{B \in \cU_n} \left|\langle 2 A + B, (2A - B) Z \rangle\right|.
    $$
    For $A \in \cU_n$ we have $\|A\|_F \le \sqrt{n}$, thus applying Cauchy-Schwarz we further get
    \begin{align*}
        L_Z(I) &\le \max_{A \in \cU_n(I)} \min_{B \in \cU_n} \left(  \| 2 A \|_F + \|B\|_F \right) \| (2 A - B) Z \|_F \\
        &\le 3 \sqrt{n} \max_{A \in \cU_n(I)} \min_{B \in \cU_n} \| (2 A - B) Z \|_F.
    \end{align*}
    Applying Lemma \ref{lemma:min_max} with linear transformations $\sigma(A) = 2AZ$ and $\tau(B) = BZ$, we obtain
    $$
        \max_{A \in \cU_n(I)} \min_{B \in \cU_n} \| (2 A - B) Z \|_F \le \max_{\substack{U \in \mathbb{R}^{n \times n} \\ \|U\|_F \le 1}} \left( \max_{A \in \cU_n(I)} \langle 2 A Z, U \rangle - \max_{B \in \cU_n} \langle B Z, U \rangle \right).
    $$
    Consider some $U$. Then
    $$
        \langle 2AZ, U \rangle = \langle 2A, U Z^\top \rangle = 2 \sum_{i \in I} \langle A_i, (UZ^\top)_i \rangle.
    $$
    Since we are estimating
    $$
        \max_{A \in \cU_n(I)} \langle 2AZ, U \rangle,
    $$
    we are free to choose every non-zero column $A_i$ independently. Note that for every $w \in \mathbb{R}^n$ we have
    $$
        \max_{\substack{u \in \mathbb{R}^n \\ \| u \|_2 \le 1}} \langle u, w \rangle = \| w \|_2,
    $$
    where the upper bound follows from Cauchy-Schwarz and the matching lower bound is achieved by setting $u_i = w_i / \|w\|$ (for $w \neq \mathbf{0}$). Thus
    $$
        \max_{A \in \cU_n(I)} \langle 2AZ, U \rangle = 2 \max_{A \in \cU_n(I)} \sum_{i \in I} \langle A_i, (UZ^\top)_i \rangle = 2 \sum_{i \in I} \|(UZ^\top)_i\|_2,
    $$
    and similarly, 
    $$
        \max_{B \in \cU_n}  \langle B Z, U \rangle = \sum_{i = 1}^n \|(UZ^\top)_i\|_2.
    $$
    Setting $\epsilon_i = 1$ if $i \in I$ and $\epsilon_i = -1$ otherwise, we conclude
    $$
        \max_{A \in \cU_n(I)} \min_{B \in \cU_n} \| (2 A - B) Z \|_F \le \max_{\substack{U \in \mathbb{R}^{n \times n} \\ \|U\|_F \le 1}} \sum_{i = 1}^n \epsilon_i \|(UZ^\top)_i\|_2.
    $$
    Note that the whole preceding discussion holds for every $I \subseteq [n]$. Our next goal is to estimate 
    $$
        \mathbb{E}_I \max_{\substack{U \in \mathbb{R}^{n \times n} \\ \|U\|_F \le 1}} \sum_{i = 1}^n \epsilon_i \|(UZ^\top)_i\|_2.
    $$
    We first turn this expression into one on which we can apply Theorem \ref{thm:ledoux}. To this end, we make use of the following elementary fact: if $g$ is a standard normal random vector in $\mathbb{R}^n$ and $v \in \mathbb{R}^n$, then
    $$
        \langle v, g \rangle \sim \mathcal{N}(0, \| v \|_2^2),
    $$
    which has the useful property
    $$
        \mathbb{E}_g |\langle v, g \rangle| = \sqrt{2/\pi} \; \| v \|_2.
    $$
    Therefore, 
    $$
        \mathbb{E}_g | \langle (UZ^\top)_i, g \rangle| = \sqrt{2/\pi} \; \| (U Z^\top)_i \|_2.
    $$
    This further implies
    \begin{align*}
        \sum_{i = 1}^n \epsilon_i \| (UZ^\top)_i \|_2 
        &= \sqrt{\pi/2}\,\mathbb{E}_g  \left[ \sum_{i = 1}^n \epsilon_i \left| (ZU^\top g)_i \right| \right] \\
        &= \sqrt{\pi/2}\,\mathbb{E}_g \left[ \|U^\top g\|_2 \sum_{i = 1}^n \epsilon_i \left| \left(Z \, \frac{U^\top g}{\|U^\top g\|_2} \right)_i \right|\right] \\
        &\le \sqrt{\pi/2}\,\mathbb{E}_g\left[ \| U^\top g\|_2 \right] \max_{\substack{h \in \mathbb{R}^n \\ \|h\|_2\le 1}} \sum_{i = 1}^n \epsilon_i \left| (Zh)_i \right|.
    \end{align*}
    Here the normalised vector is interpreted as zero when $U^\top g=0$. By Cauchy-Schwarz, we have
    $$
        \mathbb{E}_g \| U^\top g\|_2 \le \left( \mathbb{E}_g \|U^\top g\|_2^2 \right)^{1/2} = \left( \mathbb{E}_g \operatorname{Tr}( U^\top g g^\top U) \right)^{1/2} = \|U\|_F,
    $$
    where the equality follows from $\mathbb{E}_g[g_j^2] = 1$ and $\mathbb{E}_g[g_j g_k] = 0$ for $j \neq k$. As $\|U\|_F \le 1$, we finally conclude
    $$
        L_Z(I) \le 3 \sqrt{n} \cdot \sqrt{\pi /2} \; \max_{\substack{h \in \mathbb{R}^n \\ \|h\|_2\le 1}} \sum_{i = 1}^n \epsilon_i \left| (Zh)_i \right|.
    $$
     
    Again, the preceding discussion holds for every $I \subseteq [n]$. We now finally make use of the fact that $I \subseteq [n]$ is chosen uniformly at random. In particular, $\epsilon$ is a vector of independent Rademacher random variables. We are thus in a position to apply Theorem \ref{thm:ledoux}:
    \begin{align*}
        \mathbb{E}_I L_Z(I)
        &\le 6\sqrt{\pi n /2} \;\mathbb{E}\max_{\substack{h\in\mathbb{R}^n\\\|h\|_2\le1}}
            \left|\sum_{i=1}^n\epsilon_i(Zh)_i\right|\\
        &=6\sqrt{\pi n /2} \;\mathbb{E}\max_{\substack{h\in\mathbb{R}^n\\\|h\|_2\le1}}
            \left|\langle Z^\top\epsilon,h\rangle\right|.
    \end{align*}
    Applying Cauchy-Schwarz, we get
    $$
        \max_{\substack{h \in \mathbb{R}^n \\ \|h\|_2\le 1}} \left| \langle Z^\top \epsilon, h \rangle \right| \le \max_{\substack{h \in \mathbb{R}^n \\ \|h\|_2\le 1}} \| Z^\top \epsilon \|_2 \; \|h\|_2 \le \| Z^\top \epsilon \|_2.
    $$
    By Cauchy-Schwarz once again,
    $$
         \mathbb{E}_I L_Z(I) \le 6\sqrt{\pi n /2} \; \left( \mathbb{E} \|Z^\top \epsilon\|_2^2 \right)^{1/2}.
    $$
    As $\mathbb{E} \epsilon_i^2 = 1$ and $\mathbb{E} \epsilon_i \epsilon_j = 0$ for $i \neq j$, we conclude
    $$
         \mathbb{E}_I L_Z(I) \le 6\sqrt{\pi n /2} \; \|Z\|_F \le 10\sqrt{n}\;\|Z\|_F.
    $$
\end{proof}

\section{Proof of Theorem \ref{thm:main_psd}}
\label{sec: final proof}

\begin{proof}[Proof of Theorem \ref{thm:main_psd}]
    Set
    $$
        d = \min_{B \in \cS_n} \| A - B \|_F.
    $$
    We first prove a slightly stronger statement for $p=2^{-m}$, with $m\ge1$, namely
    $$
        \Pr[A_I\succeq0]\le\frac{29\sqrt{n}}{d\sqrt{p}}.
    $$

    Consider some integer $m\ge1$ and let $p=2^{-m}$. We shall sample $I\sim[n]_p$ in $m$ steps: starting with $I_0:=[n]$, in the $i$th step choose a subset $I_i\subseteq I_{i-1}$ uniformly at random. Note that the marginal distribution of each $I_i$ is $[n]_{2^{-i}}$, thus we set $I:=I_m$.

    By Theorem \ref{thm:hyperplane}, there exists a symmetric $Z\in\mathbb{R}^{n\times n}$ with $\|Z\|_F=1$ such that
    $$
        \langle A-B,Z\rangle\ge d
        \quad\text{for every }B\in\cS_n.
    $$
    Indeed, the matrix $Z$ given by Theorem \ref{thm:hyperplane} is a normalised difference of $A$ and a matrix in $\cS_n$, so it is symmetric. Let $Z'=Z-D(Z)$, where $D(Z)$ is the diagonal part of $Z$. Since $|A_{ii}-B_{ii}|\le2$, Cauchy-Schwarz gives
    $$
        \langle A-B,Z'\rangle\ge d-2\sqrt{n}
        \quad\text{for every }B\in\cS_n.
    $$
    We also have $\|Z'\|_F \le \|Z\|_F \le 1$.
    
    Following the proof outline, we first consider the sum
    \[  \sum_{j = 0}^{m-1} 4^j \left|\langle A_{I_j}-4A_{I_{j+1}},Z'\rangle\right|. \]
    Observe that we may replace the diagonal entries of $A_{I_j}$ by zeroes, since $Z'$ has zero diagonal. Thus, by Lemma \ref{lemma:subset}, for each term we have 
    \[ \mathbb{E} \left[ \left|\langle A_{I_j}-4A_{I_{j+1}},Z'\rangle\right|  \,\Big|\,  I_j \right] \leq 2 \sqrt{|I_j|}\;\|Z'_{I_j}\|_F.  \]

    If $A_I\succeq0$, that is, $A_{I_m} \in \cS_n(I_m)$, set $B^{(m)}=A_{I_m}$ and, for $j=m-1,\ldots,0$, choose
    $$
        B^{(j)}\in\operatorname*{argmin}_{Y\in\cS_n(I_j)}
            \left|\langle 4B^{(j+1)} -Y,Z'\rangle\right|.
    $$
    These minima are attained by compactness. Otherwise, if $A_I \not \succeq 0$, set $B^{(j)} = 0$ for every $j$. In either case, the triangle-inequality yields
    $$
        \langle A - B^{(0)}, Z' \rangle \, \ind\{A_I\succeq0\}
        \le\sum_{j=0}^{m-1}4^j \left|\langle A_{I_j}-4A_{I_{j+1}},Z'\rangle\right| + \sum_{j=0}^{m-1}4^j \left|\langle 4B^{(j+1)} - B^{(j)},Z'\rangle\right|.
    $$

    Since $B^{(j+1)}\in\cS_n(I_{j+1})$, we have
    $$
        \left|\langle 4B^{(j+1)} - B^{(j)},Z'\rangle\right|
        \le L_{Z'_{I_j}}(I_{j+1}),
    $$
    and so Lemma \ref{lemma:dense_psd} gives
    \[ \mathbb{E} \left[ \bigl|\langle 4B^{(j+1)} - B^{(j)},Z'\rangle\bigr|  \,\Big|\,  I_j \right] \leq 10 \sqrt{|I_j|}\;\|Z'_{I_j}\|_F.  \]
    By $\langle A - B^{(0)}, Z' \rangle \ge d - 2 \sqrt{n}$ and linearity of expectation, we obtain
    $$
        (d-2\sqrt{n})\Pr[A_I\succeq0] \le \sum_{j=0}^{m-1}4^j\mathbb{E}\Delta_j,
    $$
    where
    \[   \Delta_j = |\langle A_{I_j} - 4 A_{I_{j+1}}, Z' \rangle| + |\langle 4 B^{(j+1)} -  B^{(j)}, Z' \rangle|.  \]
    Since $\mathbb{E}|I_j|=2^{-j}n$ and
    $$
        \mathbb{E}\|Z'_{I_j}\|_F^2 = \mathbb{E} \sum_{k \neq \ell} |Z'_{k\ell}|^2 \; \ind\{k \in I_j\} \; \ind\{\ell \in I_j\} =2^{-2j}\|Z'\|_F^2\le2^{-2j},
    $$
    by Cauchy-Schwarz we have
    $$
        \mathbb{E}\Delta_j
        \le12\cdot2^{-3j/2}\sqrt{n}.
    $$
    Therefore,
    $$
        (d-2\sqrt{n})\Pr[A_I\succeq0]
        \le12\sqrt{n}\sum_{j=0}^{m-1}2^{j/2}
        =\frac{12}{\sqrt2-1}(p^{-1/2}-1)\sqrt{n}.
    $$
    Since $\Pr[A_I\succeq0]\le1$ and $2<12/(\sqrt2-1)<29$, it follows that
    $$
        d\Pr[A_I\succeq0]
        \le2\sqrt{n}+\frac{12}{\sqrt2-1}(p^{-1/2}-1)\sqrt{n}
        \le29\sqrt{n/p}.
    $$

    Finally, let $p\in(0,1]$ be arbitrary. If $p=1$, then $A_I=A$ is not positive semidefinite, so the conclusion is immediate. Otherwise, choose an integer $m\ge1$ such that $q=2^{-m}\le p<2q$. First sample $I\sim[n]_p$, and then retain each element of $I$ independently with probability $q/p$, obtaining $J\sim[n]_q$. Since $J\subseteq I$, positive semidefiniteness of $A_I$ implies positive semidefiniteness of $A_J$. Consequently,
    $$
        \Pr[A_I\succeq0]
        \le\Pr[A_J\succeq0]
        \le\frac{29\sqrt{n}}{d\sqrt q}
        \le\frac{42\sqrt{n}}{d\sqrt p}.
    $$
\end{proof}

\bibliographystyle{abbrv}
\bibliography{references}

\end{document}